\documentclass{amsart}
\def \sn {\mathbb{S}^n}

\def \crit {2^\star}

\newtheorem{lem}{Lemma}
\newtheorem{thm}{Theorem}
\newtheorem{prop}{Proposition}
\newtheorem{rmk}{Remark}

\title[Local Nirenberg problem for $Q$-curvatures]{On the local Nirenberg problem for the $Q$-curvatures}
\date{}
\author{Ph. Delano\" e, F. Robert}
\thanks{The first author is supported by the CNRS}
\begin{document}
\begin{abstract}
The local image of each conformal $Q$-curvature
operator on the sphere admits no scalar constraint although identities of Kazdan--Warner type hold for its graph.
\end{abstract}

\maketitle
\thispagestyle{empty}
\section{Introduction}
\noindent Let us call admissible any couple of positive integers $(m,n)$ such that $n>1$, and $n\geq 2m$ in case $n$ is even.
Given such a couple $(m,n)$, we will work on the standard $n$-sphere $(\sn, g_0)$ with pointwise conformal
metrics\footnote{all objects will be taken smooth}
$g_u=e^{2u}g_0$ and discuss the structure near $u=0$ of the image of the conformal $2m$-th order
$Q$-curvature increment operator
$u\mapsto \textbf{Q}_{m,n}[u]=Q_{m,n}(g_u)-Q_{m,n}(g_0)$ (see section 2), thus considering a local Nirenberg-type problem
(Nirenberg's one was for
$m=1$, \textit{cf. e.g.} \cite{mos,kw1,kw2} or \cite[p.122]{aub}). At the infinitesimal level, the situation looks as
follows (dropping henceforth the subscript
$(m,n)$):

\begin{lem}\label{lem:1} Let $L=d\mbox{{\rm \textbf{Q}}}[0]$ stand for the linearization at $u=0$ of the conformal
$Q$-curvature increment operator and
$\Lambda_1$, for the $(n+1)$-space of first spherical harmonics on $(\sn,g_0)$. Then $L$ is self-adjoint and $\hbox{Ker
}L=\Lambda_1$.
\end{lem}

\noindent Besides, the graph
$\displaystyle
\Gamma(\textbf{Q}):=\left\{(u,\textbf{Q}[u]), u\in C^\infty(\sn)
\right\}$ of \textbf{Q} in $C^\infty(\sn)\times C^\infty(\sn)$ admits
scalar constraints which are the analogue for \textbf{Q} of the so-called Kazdan--Warner identities for the conformal scalar
curvature (\textit{i.e.} when $m=1$)
\cite{kw1,kw2,bouez}. Here, a scalar constraint means a real-valued submersion defined near $\Gamma(\textbf{Q})$
in $C^\infty(\sn)\times C^\infty(\sn)$ and vanishing on $\Gamma(\textbf{Q})$. Specifically, we have:
\begin{thm}\label{kwid}
For each $(u,q)\in C^\infty(\sn)\times
C^\infty(\sn)$ and each conformal Killing vector field $X$ on $(\sn,g_0)$:
$$(u,q)\in \Gamma(\mbox{\rm {\textbf{Q}}}) \implies \int_{\sn}(X\cdot q) \ d\mu_u =0$$
where $d\mu_u=e^{nu}d\mu_0$ stands for the Lebesgue measure of the metric $g_u$. In particular, there is no solution $u\in C^{\infty}(\sn)$ to the equation:
$$Q(g_u)=z+{\rm constant}$$
with $z\in \Lambda_1$.
\end{thm}
\noindent Due to the naturality of
$Q$ (\textit{cf.} Remark
\ref{moregenatural}) and the self-adjointness of $d\textbf{Q}[u]$ in
$L^2(M_n,d\mu_u)$ (\textit{cf.} Remarks \ref{selfa} and \ref{sadj}), this theorem holds as a particular case of a general
result (Theorem \ref{jpbkw} below).
\vspace{.5em}

\noindent Can one do better than Theorem \ref{kwid}, drop the $u$ variable occuring in the constraints and find
constraints bearing on the sole image of the operator
\textbf{Q} ? Since
$L$ is self-adjoint in $L^2(\sn, g_0)$ \cite{gz}, Lemma 1 shows that the map $u\mapsto \textbf{Q}[u]$ misses infinitesimally
at
$u=0$ a vector space of dimension $(n+1)$. How does this translate at the local level? Calling now a real valued map $K$, a
scalar constraint for the local image of \textbf{Q} near 0, if $K$ is a submersion defined near $0$ in $C^\infty(\sn)$
such that
$K\circ
\textbf{Q}= 0$ near 0 in $C^\infty(\sn)$, a spherical symmetry argument (as in \cite[Corollary 5]{del1}) shows that if
the local image of \textbf{Q} admits a scalar constraint near $0$, it must admit $(n+1)$ independent such ones, that is the
maximal expectable number. In this context, our main result is quite in contrast with Theorem \ref{kwid},
namely:

\begin{thm}\label{noscalcons} The local image of {\rm \textbf{Q}} near 0 admits no scalar
constraint.
\end{thm}

\noindent Finally, the picture about the local image of the $Q$-curvature increment operator on $(\sn,g_0)$ may be completed
with a remark:
\begin{rmk}
\label{sym}
{\rm The local Nirenberg problem for \textbf{Q} near 0 is governed by the nonlinear Fredholm
formula (\ref{lfr}) (\textit{cf. infra}). In particular, as in
\cite[Corollary 5]{del1}, a local result of Moser type
\cite{mos} holds. Specifically, if $f\in C^\infty(\sn)$ is close enough to zero and invariant under a nontrivial
group of isometries of
$(\sn,g_0)$ acting without fixed points\footnote{which is more general than a free action}, then ${\mathcal D}(f)=0$ in
\eqref{lfr}, hence
$f$ lies in the local image of
\textbf{Q}.}
\end{rmk}

\medskip\noindent The outline of the paper is as follows. We first present (section 2) an independent account on general
Kazdan--Warner type identities, implying Theorem \ref{kwid}. Then we focus on Theorem \ref{noscalcons}: we recall basic
facts for the
$Q$-curvature operators on spheres (section 3), then sketch the proof of Theorem \ref{noscalcons} (section 4) relying on
\cite{del1}, reducing it to Lemma 1 and another key-lemma; we then carry out the proofs of the lemmas (sections 4 and 5),
defering to Appendice A some eigenvalues calculations. 
\section{General identities of Kazdan--Warner type}

\noindent The following statement is essentially due to
Jean--Pierre Bourguignon
\cite{bou}:
\begin{thm}
\label{jpbkw}
Let $M_n$ be a compact $n$-manifold and $g\mapsto D(g)\in C^{\infty}(M)$ be a scalar natural\footnote{in the sense of
\cite{stre}, see (\ref{topqnat}) below} differential operator defined on the open cone of riemannian metrics on $M_n$. Given
a conformal class \textbf{c} and a riemannian metric
$g_0\in \textbf{c}$, sticking to the notation $g_u=e^{2u}g_0$ for $u\in C^{\infty}(M)$, consider the operator
$u\mapsto \mbox{\rm {\textbf{D}}}[u]:=D(g_u)$ and its linearization $L_u=d\mbox{\rm {\textbf{D}}}[u]$ at $u$. Assume that,
for each $u\in C^{\infty}(M)$, the linear differential operator
$L_u$ is formally self-adjoint in $L^2(M,d\mu_u)$, where $d\mu_u=e^{nu}d\mu_0$ stands for the Lebesgue
measure of $g_u$. Then, for any conformal Killing vector field
$X$ on
$(M_n, \textbf{c})$ and any $u\in C^{\infty}(M)$, the following identity holds:
$$\int_{M}X\cdot \mbox{\rm {\textbf{D}}}[u]\, d\mu_u=0\ .$$
In particular, if $(M_n,\textbf{c})$ is equal to
$\sn$ equipped with its standard conformal class, there is no solution $u\in C^{\infty}(\sn)$ to the equation:
$$\mbox{\rm {\textbf{D}}}[u]=z+{\rm constant}$$
with $z\in \Lambda_1$ (a first spherical harmonic).
\end{thm}

\begin{proof}
We rely on Bourguignon's functional integral invariants approach and follow the proof of
\cite[Proposition 3]{bou} (using freely notations from \cite[p.101]{bou}), presenting its functional geometric
framework with some care. We consider the affine Fr\'echet manifold $\Gamma$ whose generic point is the volume form (possibly of odd type in case $M$ is not
orientable \cite{deRh}) of a riemannian metric $g\in \textbf{c}$; we denote by $\omega_g$ the volume form of
a metric $g$ (recall the tensor $\omega_g$ is natural \cite[Definition 2.1]{stre}). The metric $g_0\in
\textbf{c}$ yields a global chart of $\Gamma$ defined by:
$$\omega_g\in \Gamma \to u:=\frac{1}{n}\log\left({\frac{d\omega_g}{d\omega_{g_0}}}\right)\in C^{\infty}(M_n)$$
(viewing volume-forms like measures and using the Radon--Nikodym derivative) in other words, such that
$\omega_g=e^{nu}\omega_{g_0}$; changes of such charts are indeed affine (and pure translations). It will be easier, though, to
avoid the use of charts on $\Gamma$, except for proving that a 1-form is closed (\textit{cf. infra}). The tangent bundle to
$\Gamma$ is trivial, equal to $T\Gamma=\Gamma\times \Omega^n(M_n)$  (setting $\Omega^k(A)$ for the $k$-forms on a manifold
$A$), and there is a canonical riemannian metric on $\Gamma$ (of Fischer type \cite{fri}) given at $\omega_g\in
\Gamma$ by:
\def\lacute{\mathopen{<}}
\def\racute{\mathopen{>}}
$$\forall (v,w)\in T_{\omega_g}\Gamma,\ \lacute v,w\racute:=\int_M \frac{dv}{d\omega_g} \frac{dw}{d\omega_g}\ \omega_g\ .$$
From Riesz theorem, a tangent
covector
$a\in T^*_{\omega_g}\Gamma$ may thus be identified with a tangent vector $a^\sharp\in \Omega^n(M_n)$ or else with the function
$\displaystyle \frac{da^\sharp}{d\omega_g}=:\rho_g(a)\in C^{\infty}(M_n)$ such that:
\begin{equation}\label{covnot}
\forall \varpi\in T_{\omega_g}\Gamma,\ a(\varpi)=\int_M \rho_g(a) \varpi\ .
\end{equation}
We also consider the Lie group $G$ of conformal maps on $(M_n,\textbf{c})$, acting on the manifold $\Gamma$ by:
$$(\varphi,\omega_g)\in G\times \Gamma \to \varphi^*\omega_g\in\Gamma$$
(indeed, we have $\varphi^*\omega_g=\omega_{\varphi^*g}$ by naturality and $\varphi\in G \Rightarrow \varphi^*g\in
\textbf{c}$). For each conformal Killing field $X$ on $(M_n,\textbf{c})$, the flow of
$X$ as a map $t\in {\mathbb R} \to \varphi_t\in G$ yields a vector field $\Bar{X}$ on $\Gamma$ defined by:
$$\omega_g\mapsto \bar{X}(\omega_g):=\frac{d}{dt}\left(\varphi_t^* \omega_g \right)_{t=0}\equiv L_X\omega_g$$
($L_X$ standing here for the Lie derivative on $M_n$). 
\def\diff{\mathop{\rm Diff}\nolimits}
In this context, regardless of any Banach completion, one may
define the (global) flow $t\in {\mathbb R} \to \bar{\varphi_t}\in \diff(\Gamma)$ of $\bar{X}$ on the Fr\'echet manifold
$\Gamma$ by setting:
$$\forall \omega_g\in \Gamma,\ \bar{\varphi_t}(\omega_g):=\varphi_t^*\omega_g\ ;$$
indeed, the latter satisfies (see \textit{e.g.} \cite[p.33]{kn}):
$$\frac{d}{dt}\left(\varphi_t^* \omega_g \right)=\varphi_t^*(L_X\omega_g)\equiv
L_X(\varphi_t^*\omega_g)=\bar{X}\left[\bar{\varphi_t}(\omega_g)\right]\ .$$
With the flow $(\bar{\varphi_t})_{t\in {\mathbb R}}$ at hand, we can define the Lie derivative $L_{\bar{X}}$
of forms on $\Gamma$ as usual, by $\displaystyle
L_{\bar{X}}a:=\frac{d}{dt}\left(\bar{\varphi_t}^* a\right)_{t=0}$. Finally, one can checks Cartan's formula for
$\Bar{X}$, namely (setting $i_{\bar{X}}$ for the interior product with $\bar{X}$):
\begin{equation}\label{cartan}
L_{\bar{X}}=i_{\bar{X}}d+di_{\bar{X}}
\end{equation}
by verifying it for a generic function $f$ on $\Gamma$ and for its exterior derivative $df$ (with $d$ defined as in
\cite{lan}).\\
Following \cite{bou}, and using our global chart $\omega_g\mapsto u$ (\textit{cf. supra}), we apply (\ref{cartan}) to the
1-form
$\sigma$ on
$\Gamma$ defined at
$\omega_g$ by the function
$\rho_g(\sigma):=\textbf{D}[u]$ (see (\ref{covnot})).
Arguing as in
\cite[p.102]{bou}, one readily verifies in the chart $u$ (and using constant local vector fields on $\Gamma$) that the 1-form 
$\sigma$ is closed due to the self-adjointness of the linearized operator $L_u$ in $L^2(M_n,d\mu_u)$; furthermore (dropping
the chart $u$), one derives at once the $G$-invariance of
$\sigma$ from the naturality of
$g\mapsto D(g)$. We thus have $d\sigma = 0$ and $L_{\bar{X}}\sigma=0$, hence $d(i_{\bar{X}}\sigma)=0$ by (\ref{cartan}). 
So the function
$i_{\bar{X}}\sigma$ is constant on $\Gamma$, in other words $\displaystyle \int_M \textbf{D}[u]\ L_X\omega_u$
is independent of
$u$, or else, integrating by parts, so is $\displaystyle\int_{M}X\cdot
\textbf{D}[u]\, d\mu_u$ (where $X\cdot $ stands for $X$ acting as a derivation on real-valued functions on $M_n$).\\
To complete the proof of the first part of Theorem \ref{jpbkw}, let us show that the integrand of the latter expression at
$u=0$, namely
$X\cdot D(g_0)$, vanishes for a suitable choice of the metric $g_0$ in the conformal class \textbf{c}. To do so, we recall the
Ferrand--Obata theorem
\cite{fer,oba} according to which, either the conformal group $G$ is compact, or if not then $(M_n,\textbf{c})$ is equal to
$\sn$ equipped with its standard conformal class. In the former case, averaging on $G$, we may pick
$g_0\in
\textbf{c}$ invariant under the action of
$G$: with $g_0$ such, so is $D(g_0)$ by naturality, hence indeed $X\cdot D(g_0)\equiv 0$. In the latter case, as
observed below (section 5.1)
$D(g_0)$ is constant on $\sn$ hence the desired result follows again.\\
Finally, the last assertion of the theorem\footnote{morally consistent with Proposition \ref{natincl} (below) and Fredholm
theorem if
$L_0$ is elliptic} follows from the first one, by taking for the vector field
$X$ the gradient of $z$ with respect to the standard metric of $\sn$, which is conformal Killing as well-known.
\end{proof}

\section{Back to $Q$-curvatures on spheres: basic facts recalled}
\subsection{The special case $n=2m$} Here we will consider the $Q$-curvature increment operator given by
$\textbf{Q}[u]=Q(g_u)-Q_0$, with
\begin{equation}\label{topqcov}
Q(g_u)=e^{-2mu}(Q_0+P_0[u])
\end{equation}
where, on $(\sn,g_0)$, $Q_0=Q(g_0)$ is equal to $Q_0=(2m-1)!$ and (see \cite{br1,bec}):
\begin{equation}\label{topanexplic}
P_0=\prod_{k=1}^m\left[\Delta_0+(m-k)(m+k-1)\right],
\end{equation}
setting henceforth $\Delta_0$ (resp. $\nabla_0$) for the positive laplacian (resp. the gradient) operator of $g_0$ ($P_0$
is the so-called Paneitz--Branson operator of the metric $g_0$).
\begin{rmk}
\label{moregenatural}
{\rm One can define \cite{br2} a Paneitz--Branson
operator $P_0$ for \textit{any} metric $g_0$ (given by a formula more general than
(\ref{topanexplic}) of course), and a
$Q$-curvature
$Q(g_0)$ transforming like
\eqref{topqcov} under the conformal change of metrics $g_u=e^{2u}g_0$. Importantly then, the map $g\mapsto Q(g)\in
C^{\infty}(\sn)$ is natural, meaning (see \textit{e.g.} \cite[Definition 2.1]{stre}) that for any diffeomeorphism
$\psi$ we have:}
\begin{equation}\label{topqnat}
\psi^* Q(g)=Q(\psi^* g).
\end{equation}
\end{rmk}
\begin{rmk}\label{selfa}
{\rm From (\ref{topqcov}) and the formal self-adjointness of $P_0$ in $L^2(\sn,d\mu_0)$ \cite[p.91]{gz}, one readily verifies
that, for each $u\in C^{\infty}(\sn)$, the linear differential operator $d\textbf{Q}[u]$ is formally self-adjoint in
$L^2(\sn,d\mu_u)$.}
\end{rmk}
\subsection{The case $n\neq 2m$} The expression of the Paneitz--Branson operator on $(\sn,g_0)$ becomes \cite[Proposition
2.2]{gn}:

\begin{equation}\label{panexplic}
P_0=\prod_{k=1}^m\left[\Delta_0+\left(\frac{n}{2}-k\right)\left(\frac{n}{2}+k-1\right)\right],
\end{equation}
while the corresponding one for the metric $g_u=e^{2u}g_0$ is given by:
\begin{equation}\label{qcov}
P_u(.)=e^{-\left(\frac{n}{2}+m\right)u}P_0\left[e^{\left(\frac{n}{2}-m\right)u} . \right],
\end{equation}
with the $Q$-curvature of $g_u$ given accordingly by $\displaystyle \left(\frac{n}{2}-m\right)Q(g_u)=P_u(1)$. The analogue of
Remark
\ref{moregenatural} still holds (now see
\cite{gjms,gz}). We will consider the (renormalized) $Q$-curvature increment operator:
$\displaystyle \textbf{Q}[u]=\left(\frac{n}{2}-m\right)\left[Q(g_u)-Q_0\right]$, now with:
\begin{equation}\label{panzero}
\left(\frac{n}{2}-m\right)Q_0=\left(\frac{n}{2}-m\right)Q(g_0)=P_0(1)=\prod_{k=0}^{2m-1}\left(k+\frac{n}{2}-m\right).
\end{equation}

\begin{rmk}\label{sadj}
{\rm Finally, we note again that the linearized operator $d\textbf{Q}[u]$ is formally self-adjoint in $L^2(\sn,d\mu_u)$.
Indeed, a straightforward calculation yields
$$d\textbf{Q}[u](v)= \left(\frac{n}{2}-m\right) P_u(v) - \left(\frac{n}{2}+m\right) P_u(1)\ v\ ,$$
and the Paneitz--Branson operator $P_u$ is known to be self-adjoint  in $L^2(\sn,d\mu_u)$ \cite[p.91]{gz}.}
\end{rmk}
\vspace{.3em}
\noindent For later use, and in all the cases for $(m,n)$, we will set $p_0$ for the degree $m$ polynomial such that
$P_0=p_0(\Delta_0)$.
\section{Proof of Theorem \ref{noscalcons}}
\noindent The case $m=1$ was settled in \cite{del1} with a proof robust enough to be followed again. For completeness,
let us recall how it goes (see \cite{del1} for details).\\
If ${\mathcal P}_1$ stands for the orthogonal projection of $L^2(\sn,g_0)$ onto $\Lambda_1$, Lemma
\ref{lem:1} and the self-adjointsess of $L$ imply \cite[Theorem 7]{del1} that the modified operator
$$u\mapsto \textbf{Q}[u]+{\mathcal P}_1u$$
is a local diffeomorphism of a neighborhood of $0$ in $C^\infty(\sn)$ onto another one: set ${\mathcal S}$ for its inverse and ${\mathcal D}={\mathcal P}_1\circ\mathcal{S}$ (defect map). Then $u={\mathcal S}f$ satisfies the local non-linear Fredholm-like equation:
\begin{equation}\label{lfr}
\textbf{Q}[u]=f-{\mathcal D}(f).
\end{equation}
Moreover \cite[Theorem 2]{del1} if a local constraint exists for \textbf{Q} at $0$, then ${\mathcal D}\circ \textbf{Q}=0$ (recalling
the above symmetry fact). Fixing $z\in\Lambda_1$, we will prove Theorem \ref{noscalcons} by showing that ${\mathcal D}\circ
\textbf{Q}[tz]\neq 0$ for small $t\in {\mathbb R}$; here is how.

\noindent On the one hand, setting
$$u_t={\mathcal S}\circ \textbf{Q}[tz]:=tu_1+t^2u_2+t^3u_3+O(t^4),$$ 
Lemma \ref{lem:1} yields $u_1=0$  and the following expansion holds (as a general fact, easily verified):
\begin{equation}\label{expgen}
\textbf{Q}[u_t]+{\mathcal P}_1u_t=t^2(L+{\mathcal P}_1)u_2+t^3(L+{\mathcal P}_1)u_3+O(t^4).
\end{equation}
On the other hand, let us consider the expansion of $\textbf{Q}[tz]$:
\begin{equation}\label{expqtz}
\textbf{Q}[tz]=t^2c_2[z]+t^3c_3[z]+O(t^4)\ ,
\end{equation}
and focus on its third order coefficient $c_3[z]$, for which we will prove:
\begin{lem}\label{lem:2}
Let $(m,n)$ be admissible, then
$$\int_{\sn}z\, c_3[z]\, d\mu_0\neq 0\ .$$
\end{lem}
\noindent Granted Lemma \ref{lem:2}, we are done: indeed, the equality
$$\textbf{Q}[u_t]+{\mathcal P}_1u_t=\textbf{Q}[tz]\ ,$$
combined with (\ref{expgen})(\ref{expqtz}), yields
$$(L+{\mathcal P}_1)u_3=c_3[z],$$
which, integrated against $z$, implies:
$$\int_{\sn}z{\mathcal P}_1 u_3\, d\mu_0\neq 0$$
(recalling $L$ is self-adjoint and $z\in\hbox{Ker }L$ by Lemma \ref{lem:1}). Therefore ${\mathcal P}_1 u_3\neq 0$, hence also
${\mathcal D}\circ \textbf{Q}[tz]\neq 0$.\\
We have thus reduced the proof of Theorem \ref{noscalcons} to those of Lemmas 1 and 2, which we now
present.

\section{Proof of Lemma \ref{lem:1}}
\subsection{Proof of the inclusion $\Lambda_1\subset \hbox{Ker }L$} We need neither ellipticity nor conformal covariance for
this inclusion to hold; the naturality (\ref{topqnat}) suffices. Let us provide a general result implying at once the one we
need, namely:
\begin{prop}\label{natincl} Let $g\mapsto D(g)$ be any scalar natural differential operator on $\sn$, defined on the open cone
of Riemannian metrics, valued in $C^\infty(\sn)$. For each $u\in C^\infty(\sn)$, set
$\mbox{\rm {\textbf{D}}}[u]=D(g_u)-D(g_0)$ and $L=d\mbox{\rm {\textbf{D}}}[0]$, where $g_u=e^{2u}g_0$. Then $\Lambda_1\subset
\hbox{Ker }L$.
\end{prop}
\begin{proof} Let us first observe that $D(g_0)$ must be constant. Indeed, for each isometry $\psi$ of $(\sn,g_0)$,
the naturality of $D$ implies $\psi^*D(g_0)\equiv D(g_0)$; so the result follows because the group of such isometries acts
transitively on $\sn$. Morally, since $g_0$ has constant curvature, this result is also expectable from the theory of
riemannian invariants (see \cite{stre} and references therein), here though, without any regularity (or polynomiality)
assumption.\\ Given an arbitrary nonzero
$z\in
\Lambda_1$, let
$S=S(z)\in\sn$ stand for its corresponding ``south pole'' (where $z(S)=-M$ is minimum) and, for each small real $t$, let
$\psi_t$ denote the conformal diffeomorphism of
$\sn$ fixing
$S$ and composed elsewhere of: $\hbox{Ster}_S$, the stereographic projection with pole $S$, the dilation
$X\in\mathbb{R}^n\mapsto e^{Mt}X\in\mathbb{R}^n$, and the inverse of $\hbox{Ster}_S$. As $t$ varies, the family $\psi_t$
satisfies :
$$\psi_0=I,\qquad \frac{d}{dt}(\psi_t)_{t=0}=-\nabla_0 z$$
and if we set $e^{2u_t}g_0=\psi_t^* g_0$ we get:
$$\frac{d}{dt}(u_t)_{t=0}\equiv z.$$
Recalling $D(g_0)$ is constant, the naturality of $D$ implies 
$$\textbf{D}[u_t]=\psi_t^* D(g_0)-D(g_0) =0;$$
in particular, differentiating this equation at $t=0$ yields $Lz=0$ hence we may conclude: $\Lambda_1\subset \hbox{Ker }L$.
\end{proof}

\subsection{Proof of the reversed inclusion $\hbox{Ker L}\subset \Lambda_1$}
To prove $\hbox{Ker L}\subset\Lambda_1$, let us argue by contradiction and assume the existence of a nonzero $v\in
\Lambda_1^{\perp}\cap \hbox{Ker L}$. If ${\mathcal B}$ is an othonormal basis of eigenfunctions of $\Delta_0$ in
$L^2(\sn,d\mu_0)$, there exists an integer $i\neq 1$ and a function $\varphi_i\in\Lambda_i\cap {\mathcal B}$ (where
$\Lambda_i$ henceforth denotes the space of $i$-th spherical harmonics) such that
$$\int_{\sn}\varphi_i v\, d\mu_0\not=0$$
(actually $i\neq 0$, due to $\int_{\sn}v\, d\mu_0=0$, obtained just by averaging $Lv=0$ on $\sn$). By the self-adjointness of $L$, we may write:
$$0=\int_{\sn}\varphi_i Lv\, d\mu_0=\int_{\sn}v L\varphi_i\, d\mu_0,$$
infer (see below):
$$0=\left[p_0(\lambda_i)-p_0(\lambda_1)\right]\int_{\sn}\varphi_i v\, d\mu_0,$$
and get the desired contradiction, because $p_0(\lambda_i)\neq p_0(\lambda_1)$ for $i\neq 1$ (cf. Appendix A). Here,
we used the following auxiliary facts, obtained by differentiating \eqref{topqcov} or \eqref{qcov} at $u=0$ in the
direction of $w\in C^\infty(\sn)$:
$$\begin{array}{lll}
n=2m & \Rightarrow & Lw=P_0(w)-n!w\\
n\neq 2m & \Rightarrow & L w=\left(\frac{n}{2}-m\right)P_0(w)-\left(\frac{n}{2}+m\right)p_0(\lambda_0)w.
\end{array}$$
From $\Lambda_1\subset \hbox{Ker }L$, we get, taking $w=z\in\Lambda_1$:
\begin{equation}\label{plun}
\begin{array}{lll}
n=2m & \Rightarrow & p_0(\lambda_1)-n!=0\\
n\neq 2m & \Rightarrow & \left(\frac{n}{2}-m\right)p_0(\lambda_1)-\left(\frac{n}{2}+m\right)p_0(\lambda_0)=0.
\end{array}
\end{equation}
Moreover, taking $w=\varphi_i\in \Lambda_i$, we then have:
$$\begin{array}{lll}
n=2m & \Rightarrow & L\varphi_i=\left[p_0(\lambda_i)-p_0(\lambda_1)\right]\varphi_i\\
n\neq 2m & \Rightarrow & L \varphi_i=\left(\frac{n}{2}-m\right)\left[p_0(\lambda_i)-p_0(\lambda_1)\right]\varphi_i.
\end{array}$$

\section{Proof of Lemma \ref{lem:2}}\label{sec:lem:2}

\subsection{Case $m=2n$} For fixed $z\in \Lambda_1$ and for $t\in {\Bbb R}$ close to 0, let us compute the third order expansion of
$\textbf{Q}[tz]$. By Lemma 1 it vanishes up to first order.
Noting the identity
$$\forall v\in \Lambda_1, \frac{\textbf{Q}[v]}{Q_0} \equiv e^{-nv}(1+nv)-1\ ,$$
we find at once:
$$\frac{\textbf{Q}[tz]}{Q_0}= -2m^2t^2z^2+\frac{8}{3}m^3t^3z^3+O(t^4)\ ,$$
in particular (with the notation of section 1)
$$c_3[z]=\frac{8}{3}m^3Q_0z^3$$
and Lemma 2 holds trivially.

\subsection{Case $m\neq 2n$}
In this case, calculations are drastically simplified by picking the nonlinear argument of
$P_0$ in $P_u(1)$, namely $w:=\exp[(\frac{n}{2}-m)u]$ (see (\ref{qcov})), as new parameter for the local image of the
conformal curvature-increment operator. Since $w$ is close to 1, we further set $w=1+v$, so the conformal factor becomes:
$$e^{2u}=(1+v)^{\frac{4}{n-2m}}$$
and the renormalized $Q$-curvature increment operator reads accordingly:
\begin{equation}
\label{covtild}
\textbf{Q}[u]\equiv \tilde{Q}[v]:= (1+v)^{1-\crit}P_0(1+v)-\left(\frac{n}{2}-m \right) Q_0
\end{equation}
where $\crit$ stands in our context for $\frac{2n}{n-2m}$ (admittedly a loose notation, customary for critical Sobolev
exponents). Of course, Lemma 1 still holds for the operator
$\tilde{Q}$ (with $\tilde{L}:=d\tilde{Q}[0]\equiv \frac{\crit}{n} L$) and proving Theorem 2 (section 4) for $\tilde{Q}$ is
equivalent to proving it for \textbf{Q}. Altogether, we may thus focus on the proof of Lemma 2 for $\tilde{Q}$ instead of
\textbf{Q}\footnote{exercise (for the frustrated reader): prove Lemma 2 directly for \textbf{Q} (it takes a few
pages)}.\par
\medskip\noindent Picking $z$ and $t$ as above, plugging $v=tz$ in (\ref{covtild}), and using (from (\ref{plun})):
$$P_0(z)=p_0(\lambda_1)z\equiv (\crit-1)\left(\frac{n}{2}-m\right)Q_0z\ ,$$
we readily calculate the expansion:
$$\frac{1}{\left(\frac{n}{2}-m\right)Q_0}\ \tilde{Q}[tz]=
-\frac{1}{2}(\crit-2)(\crit-1)\ t^2z^2+\frac{1}{3}(\crit-2)(\crit-1)\crit\ t^3z^3+O(t^4)$$
thus find for its third order coefficient:
$$\frac{1}{\left(\frac{n}{2}-m\right)Q_0}\ \tilde{c}_3[z]=\frac{1}{3}(\crit-2)(\crit-1)\crit\ z^3\ .$$
So Lemma 2 obviously holds.
\appendix
\section{Eigenvalues calculations}
\noindent As well known (see e.g.
\cite{bgm}), for each $i\in\mathbb{N}$, the $i$-th eigenvalue of $\Delta_0$ on $\sn$ is equal to
$\lambda_i=i(i+n-1)$.
Recalling \eqref{panexplic}, we have to calculate
$$p_0(\lambda_i)=\prod_{k=1}^m\left[\lambda_i+\left(\frac{n}{2}-k\right)\left(\frac{n}{2}+k-1\right)\right].$$
Setting provisionally
$$r=\frac{n-1}{2}\;,\; s_k=k-\frac{1}{2},$$
so that:
$$\frac{n}{2}-k=r-s_k\, ,\; \frac{n}{2}+k-1=r+s_k\, ,\, \lambda_i=i^2+2ir,$$
we can rewrite:
\begin{eqnarray*}
p_0(\lambda_i)&=&\prod_{k=1}^m\left[(i+r)^2-s_k^2\right]\\
&=&\prod_{k=1}^m\left(\frac{1}{2}+i+r-k\right)\left(\frac{1}{2}+i+r+k-1\right)\\
&\equiv &\prod_{k=0}^{2m-1}\left(\frac{1}{2}+i+r-m+k\right),
\end{eqnarray*}
getting (back to $m$, $n$ and $k$ only)
$$p_0(\lambda_i)=\prod_{k=0}^{2m-1}\left(i+\frac{n}{2}-m+k\right).$$
In particular, we have:
$$P_0(1)\equiv p_0(\lambda_0)=\left(\frac{n}{2}-m\right) \prod_{k=1}^{2m-1}\left(\frac{n}{2}-m+k\right)$$
as asserted in (\ref{panzero}) (and consistently there with the value of $Q_0$ in case $n=2m$).
An easy induction argument yields:
$$\forall
i\in\mathbb{N},\;
p_0(\lambda_{i+1})=\frac{\left(\frac{n}{2}+m+i\right)}{\left(\frac{n}{2}-m+i\right)}\ p_0(\lambda_i)$$
(consistently when $i=0$ with \eqref{plun}), which implies: $\forall i\in\mathbb{N},\vert p_0(\lambda_{i+1})\vert > \vert
p_0(\lambda_{i})\vert$, hence in particular $p_0(\lambda_i)\neq p_0(\lambda_1)$ for $i>1$ as required in the proof of Lemma
1. Moreover, it readily implies the final formula:
$$\forall i\geq
1,\
p_0(\lambda_i)=\frac{\left(\frac{n}{2}+m\right)\ldots
\left(\frac{n}{2}+m+i-1\right)}{\left(\frac{n}{2}-m\right)\ldots \left(\frac{n}{2}-m+i-1\right)}\ p_0(\lambda_0)\ .$$

\medskip\noindent{\bf Acknowledgement:} It is a pleasure to thank Colin Guillarmou for providing us with the explicit
formula of the conformally covariant powers of the laplacian on the sphere \cite[Proposition 2.2]{gn}.\\

\vspace{.5cm}
{\begin{flushright}{Authors common address:\\\vspace{.2cm}Universit\'e de Nice--Sophia
Antipolis\\Laboratoire J.--A. Dieudonn\'e, Parc Valrose\\F-06108 Nice CEDEX
2\\\vspace{.2cm}first author's e-mail: \textsf{delphi@math.unice.fr}\\second author's e-mail: \textsf{frobert@math.unice.fr}}
\end{flushright}}

\end{document}